# Computing Pythagorean Triples

-Revised-

James M. Parks
Dept. of Mathematics
SUNY POTSDAM
Potsdam, NY USA
parksjm@potsdam.edu
July 14, 2021

***Abstract.*** Pythagorean triples are the positive integer solutions to the Pythagoras equation, $a^2 + b^2 = c^2$. They have been studied for many years, many centuries in fact. In this short paper we present a method for computing Pythagorean Triples in general, the first two cases of which go back at least to the early Pythagoreans (570-495 BCE), and to Plato (429-347 BCE). But the method may not have been investigated as extensively as presented here. Several other approaches are also considered.
This work lives at the intersection of Arithmetic, HS Algebra, and Analytical Geometry. *Dynamic Geometry* (*Sketchpad or GeoGebra)* and *Mathematica* were used to build the graphs. The work is interactive, and the reader is invited to fill in some of the details. Hints and answers are given in the final pages.

***Introduction.***

Recall the Pythagoras Theorem and its converse. While the theorem is attributed to Pythagoras, today's math researchers indicate this is not true, as it was known far before that time, for example by the Babylonians (c.1800 BCE) [8], [9].

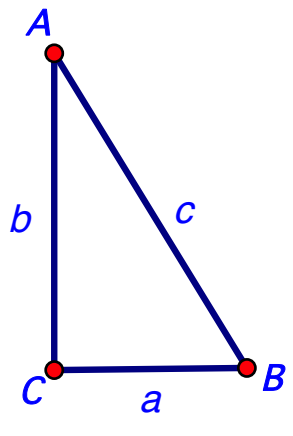

***Pythagoras Theorem.***

*Given a triangle ΔABC, with sides of length a, b, and c, such that $a \leq b < c$, the triangle is a right triangle if and only if the equation $a^2 + b^2 = c^2$ is satisfied*.

For a visual proof of the "only if" half of this Theorem, study the figure shown here. The area of each of the *2* large *(a+b)x(a+b)* squares is $(a+b)^2$. The area of the blue figure on the LH square is $a^2+b^2$. The area of the red square on the RH square is $c^2$*.* The area of the white part in each large square is *2ab.* So it can then argued that the blue area equals the red area (give an argument), and the Pythagorean equation $a^2 + b^2 = c^2$ follows [3], [4], [9].

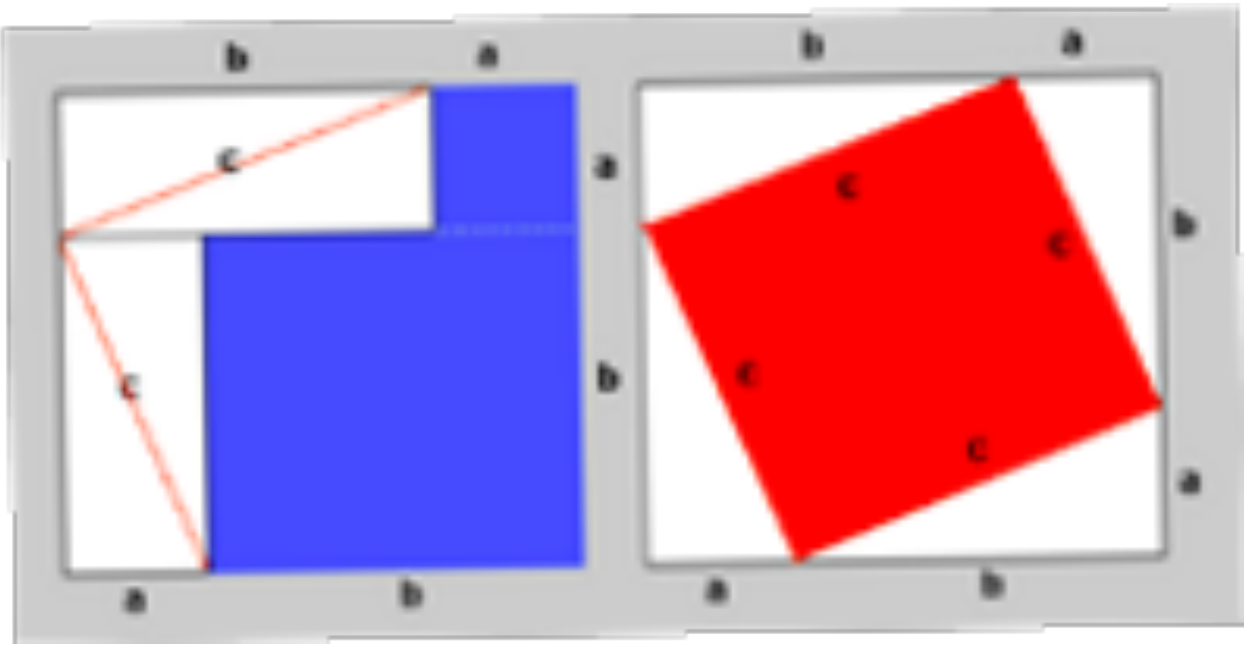

To prove the converse assume a given triangle *ΔABC,* with sides of length *a, b,* and *c,* satisfies the equation $a^2 + b^2 = c^2$. Construct a right triangle *ΔABD* with legs of length *a, b,* see the figure. Let $t = AD$, and argue that $c = t$, so the triangles are congruent by *SSS*, and thus the $<ACB$ is right (give the details).

According to Proclus (410-485 AD), who wrote a famous commentary on Euclid's *Elements*, it was Euclid (c.450-350 BCE) who gave the first complete proof of the Pythagoras Theorem, and its converse [8].

***Pythagorean Triangles.***
We are interested in studying right triangles, called *Pythagorean Triangles,* where the lengths of the sides *a, b*, *c,* are all positive integers.
For example, the well known right triangle, which has sides of lengths *3, 4, 5,* and hence satisfies the Pythagorean equation, $3^2 + 4^2 = 5^2$.

A triple of *3* positive integers $(a,b,c)$, which satisfy the Pythagoras equation, $a^2 + b^2 = c^2$, is called a *Pythagorean Triple* (PT for short). The triples *(3,4,5),* and *(6,8,10)* are PTs.

We are interested in studying Pythagorean Triples which have no common factors among the sides, that is one triple is not a multiple of another triple. Such triples are called *Primitive Pythagorean Triples* (PPT for short).
For example *(3,4,5)* is a PPT, but *(6,8,10)* is not a PPT, since *(6,8,10) = 2x (3,4,5),* see the figure.

It is customary, when stating a PT, to list the legs *a, b* first, with $a<b$, and the hypotenuse *c* last, thus we state $(a, b, c)$, although sometimes the legs may appear reversed in order, *(b, a, c),* as the result of a particular argument.

Observe that it is not possible for a PT to have $a = b$ (explain). Also, it is not possible to have $a = 1$, for then we would have to solve the equation $1+b^2 = c^2$, or equivalently $1 = c^2 - b^2 = (c - b)(c + b)$. But this is impossible (explain). A similar problem exists when $a = 2$ (why?). So the smallest value we can have for side *a* in a PT is *3,* which we know exists because we are already aware of the PPT *(3,4,5).*

Here is a list of the first *18* PPTs, ordered in size by the short leg *a*.

| ***3, 4, 5*** | ***5, 12, 13*** | ***7, 24, 25*** | *8, 15, 17* | ***9, 40, 41*** | ***11, 60, 61*** |
|---|---|---|---|---|---|
| *12, 35, 37* | ***13, 84, 85*** | ***15, 112, 113*** | *16, 63, 65* | ***17, 144, 145*** | ***19,180, 181*** |
| *20, 21, 29* | *20, 99,101* | ***21, 220, 221*** | ***23, 264, 265*** | *24, 143, 145* | ***25, 312, 313*** |

Notice that *12* of the *18* triples have the pattern *(a, b, b+1)* (listed in **bold**). Also, sides *a* and *c* are odd integers in each of these *12* PPTs, and every odd number from *3* to *25* is a value for side *a*.
This raises some questions:

- How can we find more examples of triples with this pattern *(a, b, b+1)*?
- Are the sides *a* and *b+1* always odd in PTs of the form *(a, b, b+1)*?
- Are all odd positive integers possible side *a* values?
- Are there other examples where the side *a* repeats as a side of different PTs (*a = 20* occurs twice in the list)?

We will study these questions in detail below, but first, we need to find out what is known about computing PTs. There is a general formula for computing PTs which is found in most textbooks [1]. It's called *Euclid's Formula*, but it was known long before Euclid by the Babylonians [8].

***Euclid's Formula.***

*Given a pair of integers h, k, where* $h > k > 0$*, if* $a = h^2 - k^2$*,* $b = 2hk$*, and* $c = h^2 + k^2$*, then (a, b, c) is a PT.*

The proof is a straightforward calculation and gives some insight as to why the method works. (Hint: Show $(h^2-k^2)^2 + (2hk)^2 = (h^2+k^2)^2$ ). The resulting triple (*a, b, c*) is always a PT, however it may or may not be a PPT [4].

The first problem with this Euclidean method is that not all PTs are generated by it. This is most unfortunate.
For example *(9,12,15)* is a PT, and *(9,12,15) = 3x(3,4,5),* but there do not exist *integers h, k* which will produce this triple using Euclid's formula. This is because if such *h, k* existed, then necessarily we would have $2hk = 12$, so $hk = 6$. Thus we must have either $h = 6$, and $k = 1$, or $h = 3$, and $k = 2$. But neither of these choices generates the PT *(9,12,15)* (explain why).
If we are allowed to use numbers which are not integers for *h, k*, then there exist two numbers *h, k* which will generate *(9,12,15),* namely $h = 2\sqrt{2}$ and $k = \sqrt{2}$ (show this).
But we are restricting ourselves to only positive integers here.

Fortunately, however, all PPTs are generated by Euclid's formula.
For example, if only one of *h, k* is an odd integer, $h > k$, and *h, k* are relatively prime (no common factors, other than *1*), then the triple $(h^2 - k^2, 2hk, h^2 + k^2)$ is a PPT, for the terms are also relatively prime (show this). We need only one term of *h, k* odd and the other term even, since if both are odd or both are even, the factors $h^2 \pm k^2$ are both even, hence the triple is not a PPT.

Next, given a PPT *(a, b, c)*, how do we find the pair of relatively prime integers *h, k,* such that they generate the PPT?

Since $a^2 + b^2 = c^2$ we can write $b^2 = (c - a)(c+a)$, so $(c+a)/b = b/(c - a)$. Let $h/k = (c+a)/b = b/(c-a)$, and solve for $c/b$ and $a/b$. This determines $c/b = (h^2+k^2)/2hk$, and $a/b = (h^2 - k^2)/2hk$.
We are assuming the fractions are in reduced lowest terms form here.
It then follows that since the fractions are equal, numerators are equal to numerators, and denominators are equal to denominators, so the pair $h, k$ are relatively prime. Thus the triple $(h^2 - k^2, 2hk, h^2+k^2)$ is a PPT [4].
Since this PPT is equivalent to the existence of the pair $h, k$, we have a way to find a pair of relatively prime integers which generate a given PPT using Euclid's Formula, and the Pythagorean equation.

For example, if we consider the PPT $(5, 12, 13)$, we can find $h, k$ as follows. From the equation $(c+a)/b = h/k$, we have $(c+a)/b = 18/12 = 3/2$, in reduced form. So we choose $h = 3$, and $k = 2$ (check to see that these values actually work).

The PPT form for Euclid's formula also has the equivalent form $(hk, (h^2 - k^2)/2, (h^2+k^2)/2)$, if we assume $h, k$ are both odd and relatively prime, and $h > k$. The assumption that $h, k$ are odd makes all terms even, so dividing by $2$ removes some of that. Then, when $h, k$ are relatively prime, there are no common dividers between them [1], [2], [9].

For example, in the graph here, let $h = 5$, and $k = 3$, then $a = 16$, $b = 30$, and $c = 34$, and the PT is $(16,30,34)$. Notice that $(16,30,34) = 2x(8,15,17)$, and $(8,15,17)$ is a PPT.

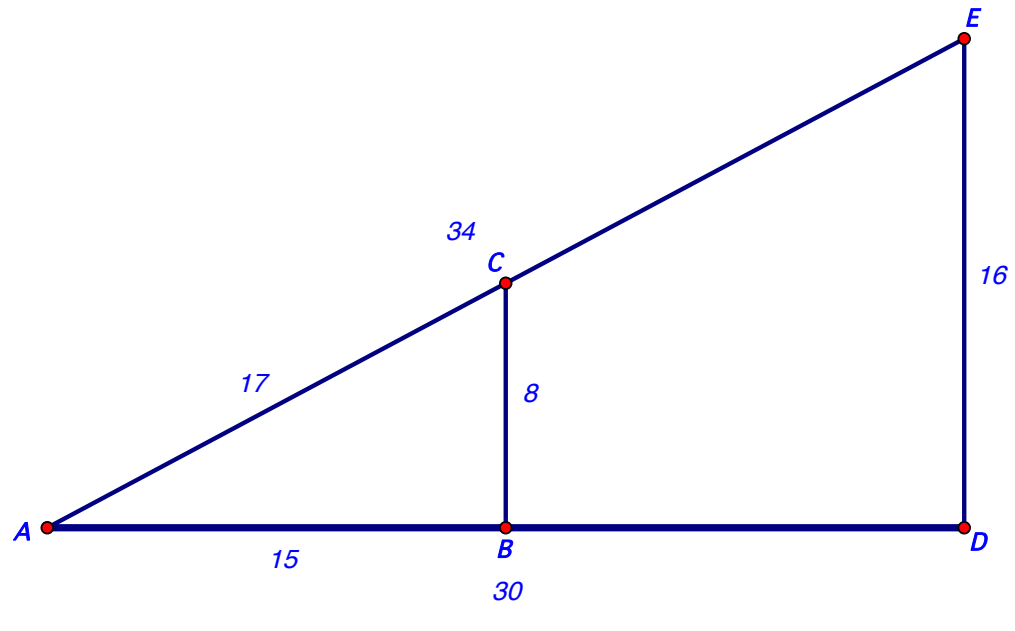

We will now investigate a method which generates triples of the form $(a, b, b+1)$, for $a$ odd. If the triple $(a, b, b+1)$ satisfies the Pythagoras equation, then $a^2 + b^2 = (b+1)^2$, so it satisfies the equivalent equation $a^2 = 2b+1$. This means $a^2$, and hence $a$, are both odd.
Solving for $b$ determines the equivalent equation $b = (a^2 - 1)/2$, which is then a positive even integer (explain).
Since all three statements are equivalent, given an odd integer $a$ you can now generate the PT $(a, b, b+1)$, by letting $b = (a^2 - 1)/2$. Thus this method works for all odd integers $a$, and we have answered at least two of the questions on p. 2.

This method is not new, in spite of the fact that it is occasional "rediscovery". In his review of Euclid's *Elements*, Proclus attributed this result to Pythagoras (the Pythagoreans formed a secret Brotherhood which credited all of their achievements to their leader Pythagoras) [8]. Proclus states that he (they) constructed the Pythagoras triple $(a, b, c)$, where $b$ and $c$ are consecutive, $c = b+1$, using the choices $a = 2n+1$, $b = ((2n+1)^2 - 1)/2$, and $c = ((2n+1)^2 - 1)/2+1$, $n$ a positive integer.
This agrees with our results above, if you let $a = 2n+1$.

For example, if $a = 3$, then $a^2 = 9$, so $b = 4$, $b+1 = 5$, and you have the PPT $(3,4,5)$. Substituting the next odd square $25$, for $a^2$ in the formula $b = (a^2 -1)/2$ gives the corresponding value $b = 12$, which then determines the PPT $(5,12,13)$.

There is a short cut calculation for these triples. If $a$ is an odd integer, compute $a^2$ and use the fact that $a^2 = b+(b+1)$ to figure out the "small half" $b$, and the "large half" $b+1$, then the triple follows.

For example, let $a = 7$, then $a^2 = 49$, so $b = 24$, and $b+1 = 25$. Hence the triple is $(7,24,25)$.

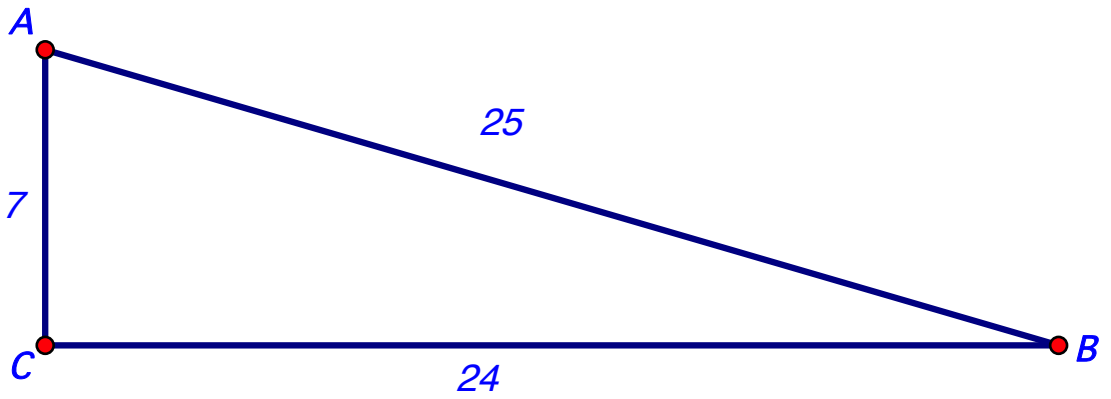

The equation choices used by the Pythagoreans above, are really Euclid's formula in disguise. If you let $h = k+1$, then $a = h^2 -k^2 = 2k+1$, $b = 2hk = 2k(k+1)$, and $c = h^2+k^2 = 2k(k+1)+1 = b+1$.
Note the method includes all primes $a = p > 2$, since primes $p > 2$ are odd.

That these PTs are actually PPTs follows from the observation that consecutive integers are always relatively prime, since they can share no common factors, except $1$ (think about it). Therefore the triple $(a, b, b+1)$ is always a PPT.
Also, note that if any $2$ numbers of a PT are relatively prime, then the PT is a PPT by the required relative prime relationship.

We will denote the *difference* between PT entries $c$ and $b$ by $d$, so $c - b = d$, and $(a,b,c) = (a, b, b+d)$. We are thus assuming $d = 1$ above.

These results are stated formally as *Theorem 1.*

***Theorem 1. PPTs with d = 1.***
*For any positive odd integer a,* $b = (a^2 -1)/2$ *if and only if* $(a, b, b+1)$ *is a PPT.*

We pause to discuss the *even/odd number* comments which appeared in the above discussions. First, it should be obvious that if you are searching for PPTs, then the $3$ integers in a triple cannot all be even, and in fact no $2$ of them can be even, for that would force the *3rd* one to be even. Likewise, the $3$ integers cannot all be odd, since the sum of any $2$ odd integers is an even integer, which would lead to a contradiction.
This leaves us with the only possibility being one entry is an even integer and two are odd integers.
The arbitrary PPT $(a, b, c)$ thus has sides $a$, $b$ which are either even, odd, or odd, even, and hypothenuse $c$ which is always odd.
That this is always the case follows since if $a$ and $b$ are both odd, say $a = 2h+1$, and $b = 2k+1$, for some $h, k > 0$. Then $a^2$ and $b^2$ are also odd, so $a^2 + b^2$ is an even integer. So $c$ must be even, $c = 2m$, for some integer $m > 0$, and $c^2 = 4m^2$. (Explain why these three assumptions lead to a contradiction.)

It would also be interesting to check for triples of the form $(a, b, b+2)$, $d = 2$, and compare these results with those we found above. So, suppose the Pythagoras equation holds for a triple of the form $(a, b, b+2)$. Then, equivalently, $b = (a^2-4)/4$, for $a$ even, and $a > 2$ (why must $a$ be even?).
For example, let $a = 4, 6, 8,$ or $10$, then $b = 3, 8, 15,$ or $24$, respectively, and these values determine the PTs ***(4, 3, 5)***, $(6, 8,10)$, ***(8, 15, 17)***, and $(10, 24, 26)$, respectively.
Notice that the *1st*, and *3rd* triples determine PPTs (in bold), where $a = 4k$, for $k = 1, 2$.
But the *2nd*, and *4th* triples determine PTs which are not PPTs, and $a = 4k + 2$, $k = 1, 2$.

The mathematical notation of denoting these two solutions for $a$ is $a = 0 \bmod 4$, in the first case, and $a = 2 \bmod 4$, in the second case.
The "*mod*" notation is short-hand from modular arithmetic. It is defined as follows. If the integer $z > 1$, and $x, y,$ are integers, then the equation $x = y \bmod z$, means $x - y = zk$, or equivalently, $x = zk + y$, for some integer $k$.
So $a = 0 \bmod 4$, means $a = 4h$, for some $h > 0$, and $a = 2 \bmod 4$, means $a = 4k+2$, for some $k > 0$. Thus the choice of $a$ as an even integer in the argument above divides the set of even integers into $2$ nonintersecting subsets, $\{4, 8, 12, \ldots\}$ and $\{2, 6, 10, \ldots\}$.
The second subset corresponds to a set of just PTs, while the first subset to a set of PPTs. For, if $a = 0 \bmod 4$, then $a = 4h$, for some $h > 0$, and $a^2 = 16h^2$, so $b = (a^2 -4)/4 = 4h^2 -1$, and $b +2 = 4h^2 +1$. The terms $b$ and $b+2$ have no common factors, since they are both odd, when $a$ is even, and they differ by $2$. Thus $(a,b,b+2) = (4h, 4h^2 -1, 4h^2 +1)$ is a PPT for all $h > 0$.
On the other hand, if $a = 2 \bmod 4$, then $a = 4k + 2$, for some $k > 0$, and $a^2 = 16k^2+16k+4$, and both are even. Then $b = 4k^2 +4k$, so $b + 2 = 4k^2 + 4k + 2$, and both of these terms are even. Ergo $(a, b, b+2) = (4k +2, 4k^2 +4k, 4k^2 +4k +2)$, for all $k > 0$.
Since $a, b,$ and $b+2$ are all even, $(a, b, b+2)$ is a PT, but not a PPT.
So, when $d = 2$, the term $a$ is always even, and if $a = 0 \bmod 4$, then $a$ determines PPTs, but if $a = 2 \bmod 4$, then $a$ determines only PTs.

Proclus attributes these results for PTs of form $(a, b, b+2)$ to Plato, a contemporary of Euclid [8]. The triple he gives has the form $a = 2n$, $b = n^2 -1$, and $c = n^2 +1$. If you let $n = 2k$, when $a = 0 \bmod 4$, and $n = 2k +1$, when $a = 2 \bmod 4$, you get our results stated above (check this).
It appears that a version of Euclid's formula can also be used to generate the PTs of the form $(a, b, b+2)$ attributed to Plato.
If you let $a = 2hk$, $b = h^2 - k^2$, and $b+2 = h^2 + k^2$, then $(h^2 + k^2) - (h^2 - k^2) = 2k^2 = 2$, so $k = 1$, and you have the PT $(2h, h^2 -1, h^2 +1)$, which is the form given by Proclus above.
If $h$ is odd, the triple is all even numbers, and $a = 2 \bmod 4$, so you have a PT which is not a PPT, as was shown above. But if $h$ is even, then $b$ and $b+2$ are odd, and you have a PPT with $a = 0 \bmod 4$, as we found above.
For example, if $h = 3$, then $a = 6$, and you have the PT $(6,8,10)$. On the other hand, if $h = 4$, then $a = 8$, and you have the PPT $(8, 15, 17)$.

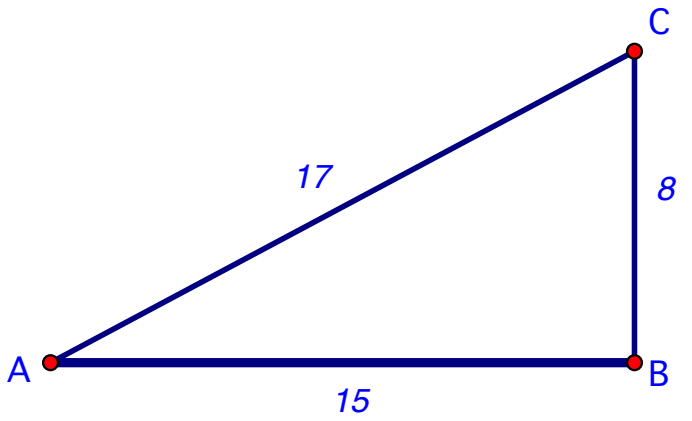

These results form *Theorem 2*.

***Theorem 2. PPTs with d = 2.***

*The triple (a, b, b+2) is a PT if and only if a is even, a > 2, and b = (a**²** - 4)/4.*
*This triple is a PPT when a = 0 mod4, but only a PT, when a = 2 mod4.*

Out of curiosity, we also check the method on the triples with form *( a, b, b+3),* where *d = 3.* Solving the corresponding Pythagoras equation for *b,* as above, you will find $b = (a^2 - 9)/6$, for *a* odd. So you must have $a = 3k$, for *k* odd, hence $a^2 = 9k^2$, $b = 3(k^2 - 1)/2$, and $b+3 = 3(k^2 + 1)/2$.
Note the similarity with the factors for *d = 1?* This is because all of the solutions for *d = 3,* are the solutions from *d = 1* multiplied by *3.*
Since all three terms have a factor of 3, all solutions are PTs, but not PPTs.

The first three triple examples for *d = 3,* are: *(9,12,15), (15,36,39),* and *(21,72,75),* which correspond to *k = 3, 5,* and *7,* respectively. Note that *(9,12,15) = 3x(3,4,5), (15,36,39) = 3x(5,12,13),* and *( 21,72,75) = 3x(7,24,25).*

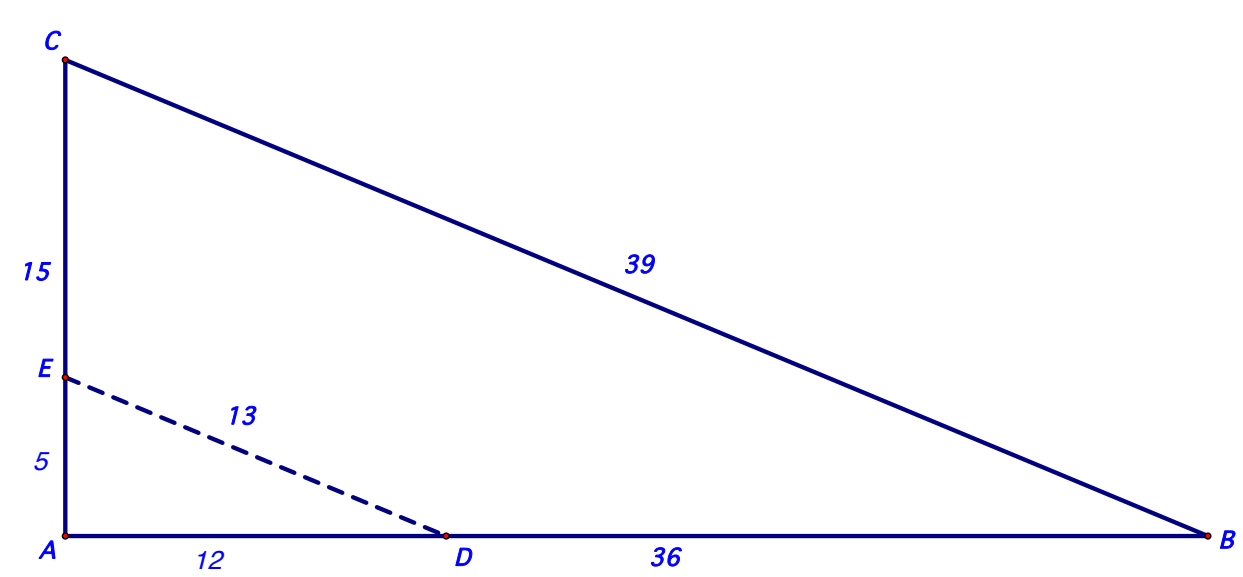

These facts makeup *Theorem 3*.

***Theorem 3. PPTs with d = 3.***

*The triple (a, b, b+3) is a PT, but never a PPT, if and only if a is odd, a = 0mod3, and b =* $(a^2 - 9)/6 = 3(k^2 - 1)/2$, *for a = 3k, k odd. Thus* $(a, b, b+3) = 3(k, (k^2 - 1)/2, (k^2 + 1)/2)$, *where* $(k, (k^2 - 1)/2, (k^2 + 1)/2)$, *k odd, is a PPT with d = 1.*

It turns out that using this general method for *d = 4, 5, 6, 7,* also determines only PTs, while the choices for *d = 8, 9* determine both PTs and PPTs in each case, similar to the results for *d = 2* above (check these claims as time allows). The observation above for *d* = 3 applies here also. The solutions *4, 5, 6, 7,* are mostly multiples of the results for *d = 1* or *2*.

You may be interested to learn that out of the first *100* PPTs, the associated *d* values for PPTs are *1, 2, 8, 9, 18, 25, 32, 49, 50, 72, 81, 98, 121, 128,* and *162,* see *A096033*, with the larger number of examples appearing for the smaller values of *d*, the majority of the examples being for *d = 1* and *2,* [4], [5].
We call the set of all values of *d* for PPTs in *A096033*, the *allowable values of d.*
Note that the odd values of *d* are the squares of the odd integers, and the even values are twice the squares of the consecutive integers.
The list of PPTs appears to become somewhat sparse in the set of all PTs as *d* increases.

Using the method of a fixed value for *d* is handy, if you only want to find PTs. The surprising fact is that there are no primitive Pythagoras triangles, PPTs, which have the difference between the longest leg and the hypotenuse equal to one of the values not on the list of allowable *d's* above.

We now consider the general construction of PTs for arbitrary values of $d$.
As we saw in the examples above, there is a calculation which provided an answer to the search for a PT when $d = 1, 2,$ and $3,$ so we will generalize this method.
If we let $(a, b, b+d)$ be a PT, with $d$ a positive integer, such that the Pythagoras equation holds, then, equivalently, $b = (a^2 - d^2)/2d,$ and $b+d = (a^2 + d^2)/2d$.
So $a$ and $d$ must both be either odd or even. Also $a^2$ must have $d$ as a factor.
Therefore, given any $d$, choose a value for $a$ which is a multiple of $d,$ say $dm,$ $m$ a positive integer which is arbitrary if $d$ is even, and odd if $d$ is odd. Then the triple $(a,b,b+d)$ will be a PT, whenever $b = (a^2 - d^2)/2d,$ and conversely.

These facts form *Theorem 4.*

***Theorem 4. PPTs with d an Arbitrary Positive Integer.***
*Suppose a, d are both odd or both even positive integers, and d divides* $a^2$*, then the triple (a, b, b+d) is a PT if and only if* $b = (a^2 - d^2)/2d$.

We now consider some interesting special cases of PPTs.

Speaking of prime numbers, you may have noticed that in many of the above examples of PPTs, with $d = 1,$ the short leg '*a'* is a prime number. After all, prime numbers greater than $2$ are odd integers.
For example, the PPTs *(3,4,5), (5,12,13),* and *(7,24,25),* all have prime numbers for the length of leg $a$.
Surprisingly, the only PPTs with short leg '*a*' a prime number are those with $d = 1.$
This follows since, if $d > 1,$ and $d$ is odd, then $d = n^2,$ $n$ odd, $n > 1,$ and by empirical evidence $a = n \bmod 2n$, so $a = (2k+1)n,$ for $k \geq 1,$ and $a$ is not prime.

It also happens occasionally that a PPT *(a,b,c)* will have prime numbers in both the $a$ leg position, and the hypotenuse position, such as *(3,4,5)* and *(5,12,13)* [7]. These PPTs are sometimes called *prime* or *2-prime* PPTs [16].

If *(a,b,c)* is any PPT, it follows that one of the legs $a, b$ is divisible by $3,$ and one of them is divisible by $4,$ while one of the sides $a, b, c,$ is divisible by $5$ [17].
Examples are *(3,4,5), (5,12,13),* and *(7,24,25).*

That one of $a, b$ is divisible by $3$ follows by assuming $a, b$ are not divisible by $3.$ Then when you divide $a, b$ by $3,$ you have $a = 3h + 1,$ or $3h + 2$, and $b = 3i + 1,$ or $3i + 2.$ These statements can be restated as $a = 3h \pm 1,$ for some integers $h > 0,$ and $b = 3i \pm 1,$ for some integers $i > 0.$
Thus $a^2 + b^2 = (9h^2 \pm 6h + 1) + (9i^2 \pm 6i + 1) = 3(3h^2 \pm 2h + 3i^2 \pm 2i) + 2.$ This is not the square of a positive integer, since it is not divisible by $3$. When you divide $c$ by $3$ you have $c = 3j \pm 1,$ for some integers $j > 0,$ and $c^2 = 3(3j^2 \pm 2j) + 1.$
Thus we have a contradiction, since the remainders $2$ and $1$ are not equal.

A similar argument holds for the divisor $5$.

We know that one number $a$, $b$ must be divisible by $4$, by the following argument.
Since we have $b^2 = c^2 - a^2 = (c - a)(c + a)$, these last two numbers must be even, for they are the difference and sum of two odd numbers, so $(c - a) = 2m$, and $(c + a) = 2n$, for positive integers $m$, $n$. Thus we have $b^2 = 4mn$, and we can solve the system for $c$ and $a$, which determines $c = m + n$, and $a = n - m$. Also, it follows that $m$ and $n$ are relatively prime (show this). We have $b$ even, $b = 2k$, for some positive integer $k$, so $4k^2 = 4mn$, and $k^2 = mn$. Thus we must have $m = s^2$, and $n = t^2$, so $m$ and $n$ are relatively prime. By the above this means that $c = s^2 + t^2$, and $a = s^2 - t^2$, then we have $b = 2st$, so now $s$, $t$ are Euclidean variables generating $(a, b, c)$, and cannot both be even, since they are relatively prime, and they can't both be odd, since then $a$ would be even, contradiction. Thus $b$ is divisible by $4$.

As a consequence of this result, if $(a,b,c)$ is a prime PPT, with $a$ and $c$ primes, $a > 7$, then $60|b$. For example $(11,60,61)$, $(19,180,181)$.
Thus if $(a,b,c)$ is a PPT with $a = p$, a prime, with $p > 7$, and we suspect $(a,b,c)$ is a prime PPT, we can check to see if $60$ divides $b$. If not, then $(a,b,c)$ is not a prime PPT.
A similar argument holds for $(a,b,c)$ with $c = p$, a prime. If $60$ does not divide $b$, then $a$ is not a prime, so $(a,b,c)$ is not prime.

Let $a = p$, $p$ a prime, $p > 3$, then $(p, b, b+1)$ a PPT, $d = 1$, and since $p$ is odd, $b = (p^2 -1)/2$, and $b+1 = (p^2 +1)/2$. By the factorization $p^2 -1 = (p +1)(p -1)$, it follows that one factor is divisible by $2$ and the other by $4$, since the $3$ numbers $(p -1)$, $p$, $(p +1)$ are consecutive integers, and $p$ is odd. Also, one factor is divisible by $3$, and it's can't be $p$. Thus $2b = p^2 -1$ is divisible by $2x3x4 = 24$, so $b$ is divisible by $12$, and $b+1 = 1mod4$.

Numbers with the form *1mod4* are called *Pythagorean Primes,* but they may not be prime, for example $9 = 1mod4$. See also the examples above, and Fermat's Theorem below [13].

If we use the Euclid formulation of the coordinates $a$, $b$, $c$, as $h^2 - k^2$, $2hk$, and $h^2 + k^2$, then in order for both $a$ and $c$ to be prime, we must have $h^2 - k^2 = (h -k)(h +k)$ prime, so then $h - k = 1$. This leads to the result that $a = 2k+1$ and $c = 2k^2+2k+1$ and both are prime.

There are also PPTs, where the hypotenuse $c$ alone is prime, such as $(9,40,41)$, $(15,112,113)$, and $(21,20,29)$. These PPTs have the form $(a, p-d, p)$, where $d \geq 1$, $p$ is prime, but $a$ is not prime. We construct a few examples.
Let '$a$' be an odd integer multiple of $5$, $a = 5(2k+1)$, $k > 0$, $k \neq 12$, and $k \neq 0mod13$. Then $(a,c -1,c)$ is a PPT with $d = 1$, where $c -1 = (a^2 -1)/2 = (25(2k+1)^2 -1)/2$, and $c = (25(2k+1)^2 +1)/2$. It follows that $c$ always ends in $13$ (check this), and $c$ is often a prime number (why are there extra restrictions on $k$?).
For example, if $a = 15, 25, 35,$ or $45$, then $(a,b,c) = (15,112,113)$, $(25,312,313)$, $(35,612,613)$, or $(45,1012,1013)$, resp., and all of these $c$ terms are prime numbers which end in $13$. Unfortunately, $c$ is not always prime here, for if $a = 55$, then $c = 17x89$.

Similarly, let the positive odd integer '*a'* end in a *1* or *9,* and consider the PPT *(a,c -1,c), d = 1*, and *c -1 = (a² -1)/2.* Then *c* always ends in a *1,* and *c* is sometimes a prime number*.*
For example, if *a = 39, 49, 51, or 121,* then the corresponding PPTs are *(39, 760, 761), (49, 1200, 1201), (51, 1300, 1301),* and *(121, 7320, 7321),* resp*.,* and *761, 1201, 1301,* and *7321* are all prime numbers (check these). However, if *a = 99* or *81,* then the generated triples are the PPTs *(99, 4900, 4901)* and *(81, 3280, 3281),* resp., but neither *4901* nor *3281* is a prime number (check this).
Unfortunately, if term *a* ends in a *3* or *7,* and *a > 3, d = 1,* then the PPTs all have terms *c* which end in a *5,* so *c* cannot be prime.

The existence of PPTs which have a leg and the hypotenuse both prime form a very scarce set. These PPTs are called *prime PPTs* [16]. The first few members are *(3,4,5), (5,12,13), (11,60,61), (19,180,181*), and *(29,420,421)*.
From there the numbers get large very quickly, see *A048161* [7], [16]. These triples all correspond to very skinny right triangles.

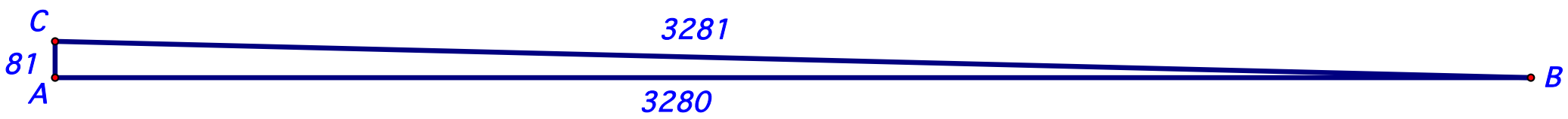

Interestingly, the distribution of primes decreases as the numbers increase in size, and the percent of prime numbers in the set of prime integers *{1, 3, …, n}* slowly approaches *0* as *n* increases. This is formally called the *Prime Number Theorem* and it states that the number of primes less than or equal to *N* is approximately *N/ln(N)*. Unfortunately the proof of this result is beyond the scope of this work [13].

Notice that the hypothenuses *c* of all of the PPTs which we have dealt with so far satisfy *c = 1mod4,* so *c = 4n+1,* for some positive integer *n.* This turns out to be true of all PPTs, but we cannot give a proof here, as it involves math outside the scope of this paper, see [9]. Integers of the form *4n+1* are called *Pythagorean primes,* but they need not be real prime numbers, as we saw in the examples above. If they are also prime numbers, they are exactly the numbers which are the sum of *2* squares of integers, by the following result of Fermat [10], [11].

***Fermat's Theorem on the Sums of Two Squares.***

*Every prime p is the sum of two squares if and only if p is of the form 4n+1, for some positive integer n*.

Thus there is a connection between Pythagorean Triples and Fermat's Theorem. If you consider a PPT in the form of Euclid's formulas, the hypotenuse is the sum of two squares, $c = h^2 + k^2$. If *c* is also prime, then we must have *c = 1mod4,* by Fermat's Theorem*.* The problem is that there also exist nonprime hypotenuses *c* which have the form *4n+1,* such as *(7,24,25).*

***Corollary, Theorem 4. PPTs with d Arbitrary.***

a.) *If (a,b,c) is a PPT, then either one of a, b is divisible by 3, and either one by 4, and one of a, b, c is divisible by 5.*

b.) *If a = p a prime number and p > 3, then the triple (p, b, b+1) is a PPT with d = 1, b = $(p^2 -1)/2$ is divisible by 12, b+1 = 1mod4, and b+1 is often a prime number.*

c.) *If a = 5(2k+1), k > 0, k≠12, and k ≠ 0mod13, then (a,c -1,c) is a PPT, where c always ends in 13, and c is often a prime number.*

d.) *If an odd positive integer a ends in 1 or 9, in the PPT (a, c -1, c) with d = 1, where c = $(a^2 +1)/2$, then c always ends in a 1, and c is often a prime number.*

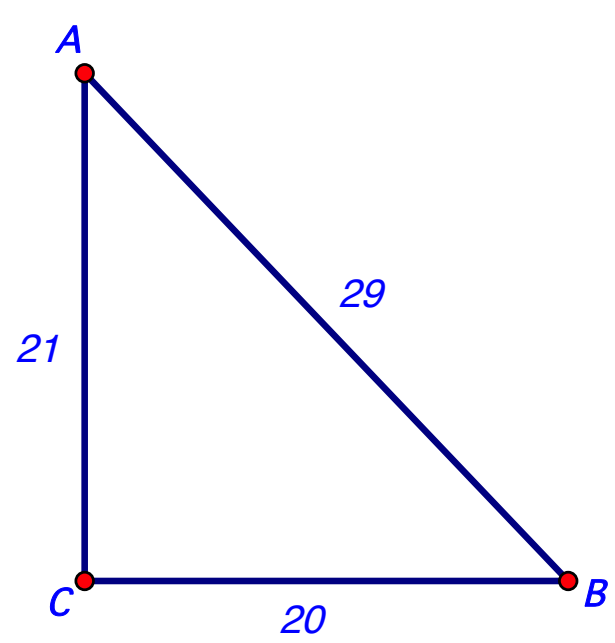

*Theorem 4* only determines solutions for PPTs where the hypotenuse and the longest leg differ by length *d*.
Consider the example shown here of a right triangle with one leg of length *20* and the other leg of length *21*. Then $20^2 + 21^2 = 29^2$ (check this), and *(20,21,29)* is a PPT (why?).
This opens up a new case of PPTs where the legs differ by *1* unit, and the PPT has the form *(a, a+1, c)*, which satisfies the equation $2a^2 + 2a + (1 - c^2) = 0$, for *c* odd.
Solve this equation for *a,* and you have $a = (\sqrt{(2c^2 -1)} -1)/2.$ So you are searching for odd integers *c* such that $2c^2 -1$ is a perfect square. This is not a very efficient approach.
An easier approach is to double the *Pythagoras* equation $a^2 + (a+1)^2 = c^2,$ which then simplifies to $(2a+1)^2 +1 = 2c^2.$ Now we are looking for odd squares $(2a+1)^2$ such that $((2a+1)^2 +1)/2,$ is a perfect square, see *A001653* [5].
Such numbers are few and far between.

A recursion formula for the sequence of the values of sides *c* in the PPT *(a, a+1, c)* is known, see *A001653* [5]. The first *4* examples of *c* are *5, 29, 169*, and *985.* The corresponding PPTs are *(3, 4, 5), (20, 21, 29), (119, 120, 169),* and *(696, 697, 985)*. The numbers get large very quickly.

Here is a better method for finding triples of the form *(a,a+1,c)*, which uses a combination of Euclid's equation, and a recursion formula [4].
For *h, k* the values in Euclid's formula, the smallest *1st* pair for this form is *(h, k) = (2, 1),* which generates the PPT *(3, 4, 5)*, and the *2nd* smallest pair is *(5, 2),* which generates the PPT *(20, 21, 29).* Then the next pair *(h',k')* has the value for *h'* obtained by multiplying the *h* value *5* for the *previous pair* by *2*, and add the *h* value *2* for the *previous pair* to that result: *h' = 2x5 + 2 = 12.* You use the same formula for *k': k' = 2x2 +1 = 5,* in this case*.* The pair *(12, 5)* generates the next PPT *(119,120,169),* and so on.

Here is yet another recursion method which generates PPTs of the form *(a, a+1, c)* [4].
Suppose you know a PPT and the previous smaller PPT, which both have the form *(a, a+1, c),* and they are ordered by the leg size of *a.* Such as *(20, 21, 29)* and *(3, 4, 5).* Then the next larger PPT can be calculated by multiplying the smallest leg size of the *last* PPT, *20,* by *6*, subtract the smallest leg size of the *previous* PPT from that, *3*, and

add *2:* $(6x20 - 3) + 2 = 119$. It follows that $c^2 = 120^2 + 119^2 = 169^2$, and the new PPT is *(119, 120, 169)*.
The hypotenuse *c* can also be calculated using the recursion formula: multiply the hypotenuse *29* of the *previous* PPT by *6* and subtract the hypotenuse *5* of the PPT *before* that: $6x29 - 5 = 169$.

These results are stated as *Theorem 5.*

***Theorem 5. Triples (a, a+1, c).***

a.) *If (a, a+1, c) is a PPT, then* $a = [-1 + \sqrt{(2c^2 -1)}]/2$, *and conversely.*

b.) *If the generating pairs (h, k) in the Euclidean formula for two consecutive triples of type (a, a+1, c), ordered by the leg a are known, say (h, k) and (h', k'), then the next triple of this type is generated by (h", k"), where* $h'' = 2h' + h$, *and* $k'' = 2k' + k$.

c.) *If 2 consecutive triples (a, a+1, c) and (a', a'+1, c'), ordered by the leg a are known, then the next triple in the sequence (a", a"+1, c") is generated by* $a'' = 6a' - a + 2$, *where c" can be computed by the Pythagoras equation, or by the recursion formula* $c'' = 6c' - c$.

By now you are aware that there are many formulas for computing PTs. Here is another one for computing PTs in general which is simple, and indirectly related to the above techniques for finding triples of the form *(a, b, b+d)* [4]. But instead of choosing *d,* the short leg *a* is fixed, and all of the possible Pythagorean Triples which share that side are computed, using the integral factors of $a^2$ to find the possible values of *d.*

Here's the method: Fix *a* as a known positive integer, solve the equation $a^2 + b^2 = c^2$ for $a^2$, then factor, to get the equation $a^2 = (c + b)(c - b)$.
Now, how can you determine *b* and *c* so that you have a PT?
If there exist *2* factors of $a^2$, say *u, v,* with $v > u$, then $a^2 = uv$, and $u < a$. Therefore, let $v = c + b$, and $u = c - b$. Solve for *b* and *c* in terms of *u and v,* to find $c = (u + v)/2$, and $b = (v - u)/2$. These values determine the PT *(a, b, c).* Notice also that *u* and *v* must both be either even or odd (why?).
For example, let $a = 20$, then $a^2 = 400$, and $400 = 2^4x5^2$. The allowable factorizations of *400* are: *2x200, 4x100, 8x50,* and *10x40* (recall the restrictions on *u & v*).
If $u = 2$, and $v = 200$, then $c = (u+v)/2 = 101$, and $b = (v-u)/2 = 99$, so the PPT is *(20, 99, 101),* and $d = 2$.
Similarly, for $u = 4$, and $v = 100$, the PT *(20,48,52),* is determined, where $d = 4$.
For $u = 8$, and $v = 50$, the triple is *(20, 21, 29),* the PPT with $d = 8$.
For $u = 10$, and $v = 40$, the PT *(20,15, 25)* is determined, and $d = 5$.
So, there are only *4* PTs with $a = 20$, and *2* of them are PPTs, see the figure.

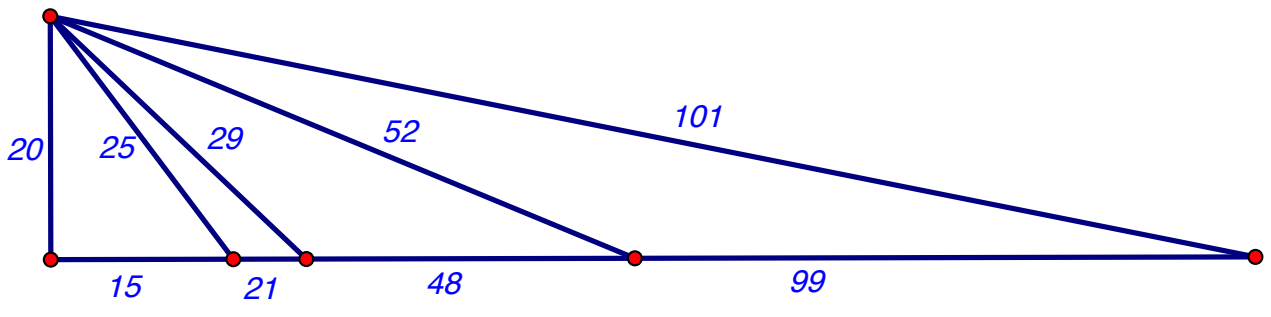

These facts are stated as *Theorem 6.*

***Theorem 6. Given a, Determine u, v.***

*In order to generate a PT (a, b, c), with side ‘a’ known, determine factors u, v, both even or both odd, and v > u, such that* $a^2 = uv$*, then let* $b = (v - u)/2$ *and* $c = (u + v)/2$*.*

A simple change in the terms of the Pythagoras equation gives the following variation on the *u, v* method which assumes both side *a* and *d* are known for the triple *(a, b, b+d)* [4]. Solve the equation $a^2 + b^2 = (b + d)^2$ for $a^2$, to obtain $a^2 = d(2b + d)$.
Then *d* must be less than *a,* so all factors of $a^2$ which are less than *a* are possible values for *d,* that is *d* must satisfy: *d* divides $a^2$, and $d < a$.
Since *a* and *d* are assumed to be known, you can solve for *b* in terms of *a* and *d,* $b = (a^2/d - d)/2$. Note $(a^2/d - d)$ must be even.

For example, let side $a = 12$, then $a^2 = 144 = 2^4 x 3^2$. The ‘possible’ values for *d* are then: *1, 2, 3, 4, 6, 8,* and *9*. Checking each *d* value in the formula for *b* we have the following results.

| *d* | $(a^2/d - d)$ | *even/odd* |
|---|---|---|
| *1:* | *144/1 - 1 = 143*, | odd. |
| *2*: | *144/2 - 2 = 70*, | even, *b = 35*, *c = 37*, and we have the PPT *(12, 35, 37).* |
| *3:* | *144/3 - 3 = 45*, | odd. |
| *4:* | *144/4 - 4 = 32,* | even*, b = 16, c = 20,* and we have the PT *(12, 16, 20).* |
| *6*: | *144/6 - 6 = 18,* | even, *b = 9, c = 15,* and we have the PT *(12, 9, 15).* |
| *8*: | *144/8 - 8 = 10,* | even*, b = 5, c = 13,* and we have the PPT *(12, 5, 13)* |
| *9*: | *144/9 - 9 = 7,* | odd. |

Hence, there are *4* PTs, and *2* are PPTs, one with $d = 1$, and one with $d = 2$.

These general facts are stated as *Theorem 7*.

***Theorem 7. Given a, d, Determine (a,b,c).***

*If a and d are known positive integers, such that* $d < a$*, and* $a^2/d - d$ *is even, then when* $b = (a^2/d - d)/2$*, the PT (a, b, c) is determined.*

One way to analyze objects with a geometric connection is to graph them. This can be done with PTs and right triangles as follows. The PT *(a,b,c),* $a<b<c$, determines a triangle which can be graphed by plotting vertices on an *xy* coordinate plane. The points *O*, *(a,0),* and *(a,b),* give the vertices of a right triangle congruent to the right triangle for *(a,b,c)* (why?), see the figure. The vertex *(a,b)* is the key to the graph of this PT, if that point is graphed on the *xy* coordinate plane by itself, then there is enough information, by the above observation, to reconstruct the triangle [4].

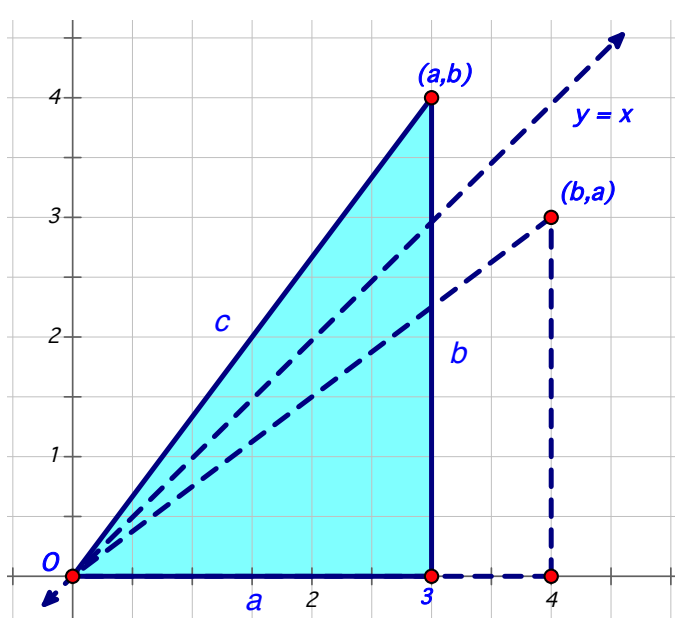

Since the PTs *(a,b,c)* and *(b,a,c)* are both legitimate representatives of the same triangle, both pairs *(a, b)* and *(b, a)* will be plotted. The plot of all such pairs of points determined by the set of all PTs are symmetric about the line $y = x$, by reflection, see the figure. The red points have $a < b$, and black points have $a' > b'$. The plot graph here is of all PTs with legs $a,b < 1000$ [4], [12], [14].

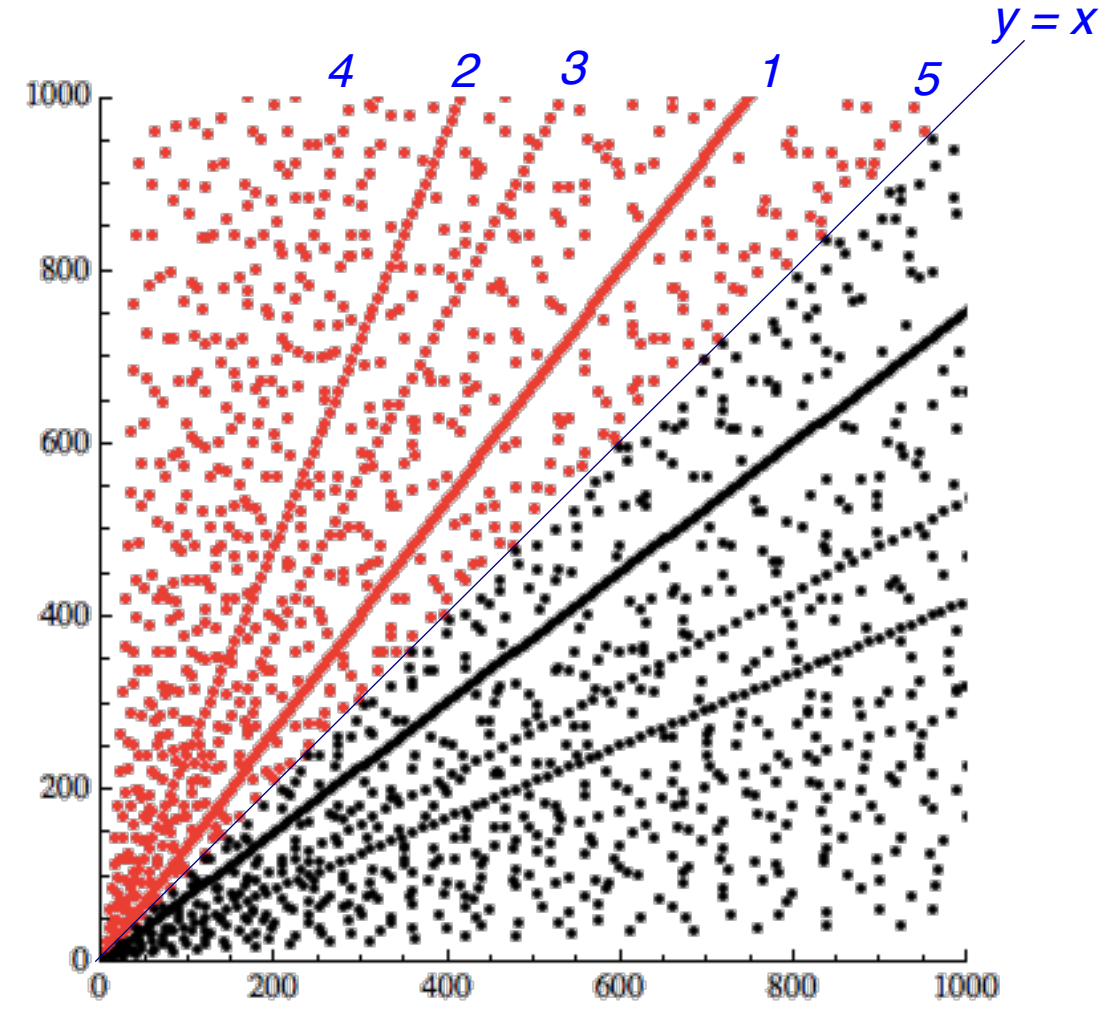

The *2* darkest straight linear arrays of points next to the line y = x are the integer multiples of the points *(3,4),* and *(4,3).* These linear arrays are more dense, since *(3,4,5)* determines the generating points *(3,4), (4,3)* which are closer to the origin *O.* The "lines" are labeled *1* on the red side of the graph, the dark black "line" is the reflection of "line *1*". The next less dense linear array, labeled *2,* corresponds to the PPT *(5,12, 13).* This is followed by the slightly less dense linear array, labeled *3*, corresponding to the PPT *(8,15,17)*. Next is the still less dense linear array, labeled *4*, which correspond to the PT *(7, 24, 25)*. The very faint linear array, labeled *5*, which is close to $y = x$, corresponds to the PPT *(20,21,29)*.

If you remove the nonPPTs from this plot you have the PPTs, see the next figure shown below [4]. There are some faint curves and patterns appearing. If you plot only the PPTs with legs $a,b < 10,000$, you get the even more striking figure below [4], [14]. The patterns are truly spectacular, and suggest hidden connections among the actual list of all PPTs. This plot graph was apparently first produced by R. Knott using *Mathematica* [4].

Some parabolic shaped curves that can be seen in this figure can be explained by the *d* factors. These parabolic curves will contain points with the first *2* coordinates of the PPTs for each allowable value of *d.*

The parabolas for $d = 1$, correspond to the PPTs for this value of *d.* These PPTs have the form $(a,(a^2-1)/2,$ and $(a^2+1)/2)$, for *a* odd, so the equations of the parabolas are given by these formulas. Can you give the equations? All you have to do is let $x = a$, and $y = b$, thus the equation is $y = (x^2 - 1)/2$.

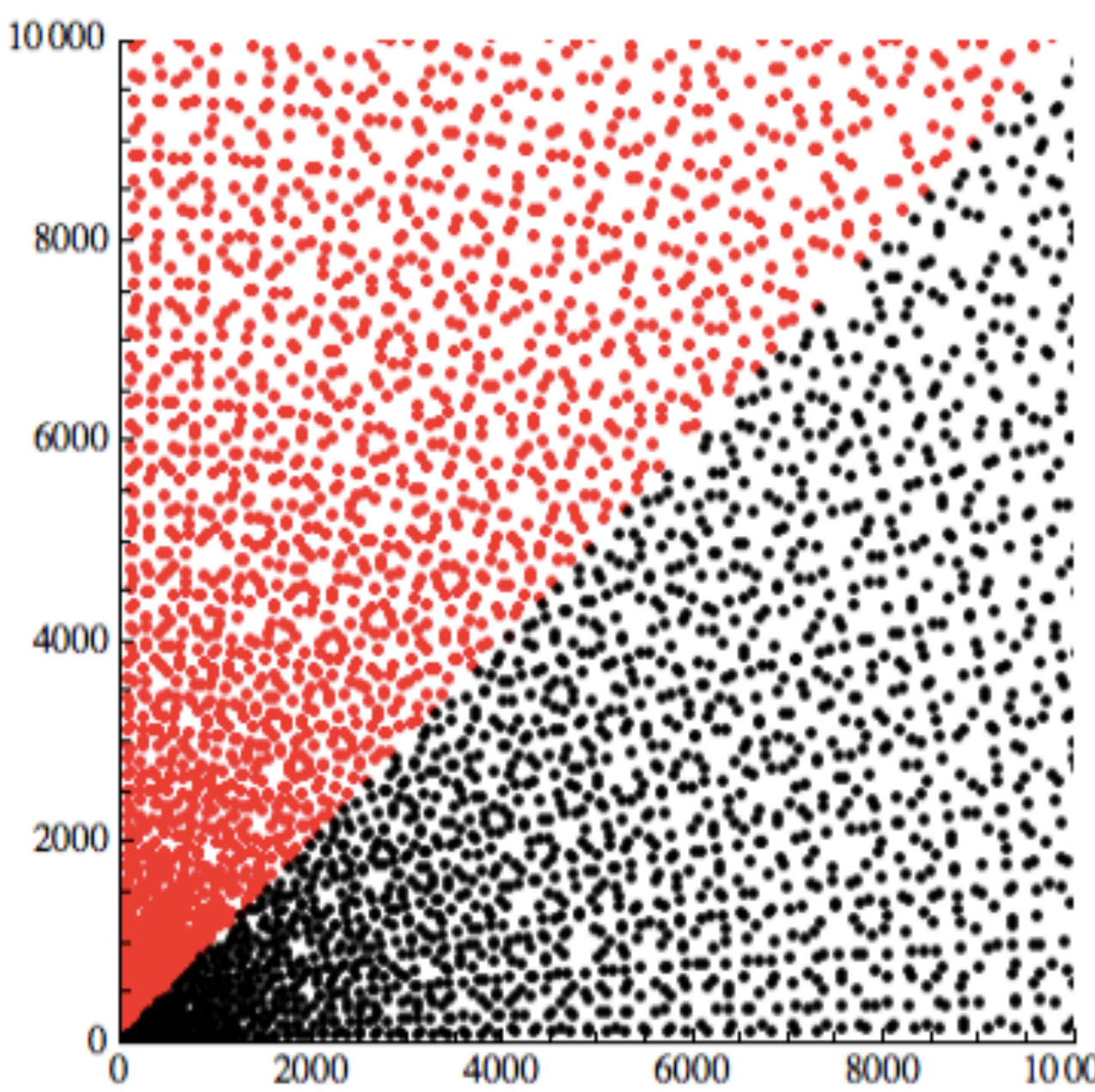

Similarly, for the reflected point, we let $y = a$, and $x = b$, so the equation is $x = (y^2-1)/2$. These parabolas are reflections of each other about the line $y = x$.

The first equation is a parabola which opens up about the positive $y$-axis, the second equation is a parabola which opens to the right about the positive $x$-axis. The points $(a,(a^2-1)/2)$ and $((a^2-1)/2,a)$, are the corresponding points of the PPT, and are on the respective parabolas.

For example, since $(3,4,5)$ is a PPT with $d = 1$, the two parabolas for $d = 1$ which contain the points $(3,4)$, $(4,3)$ are $y = (x^2-1)/2$, and $x = (y^2-1)/2$, respectively, see the graph.

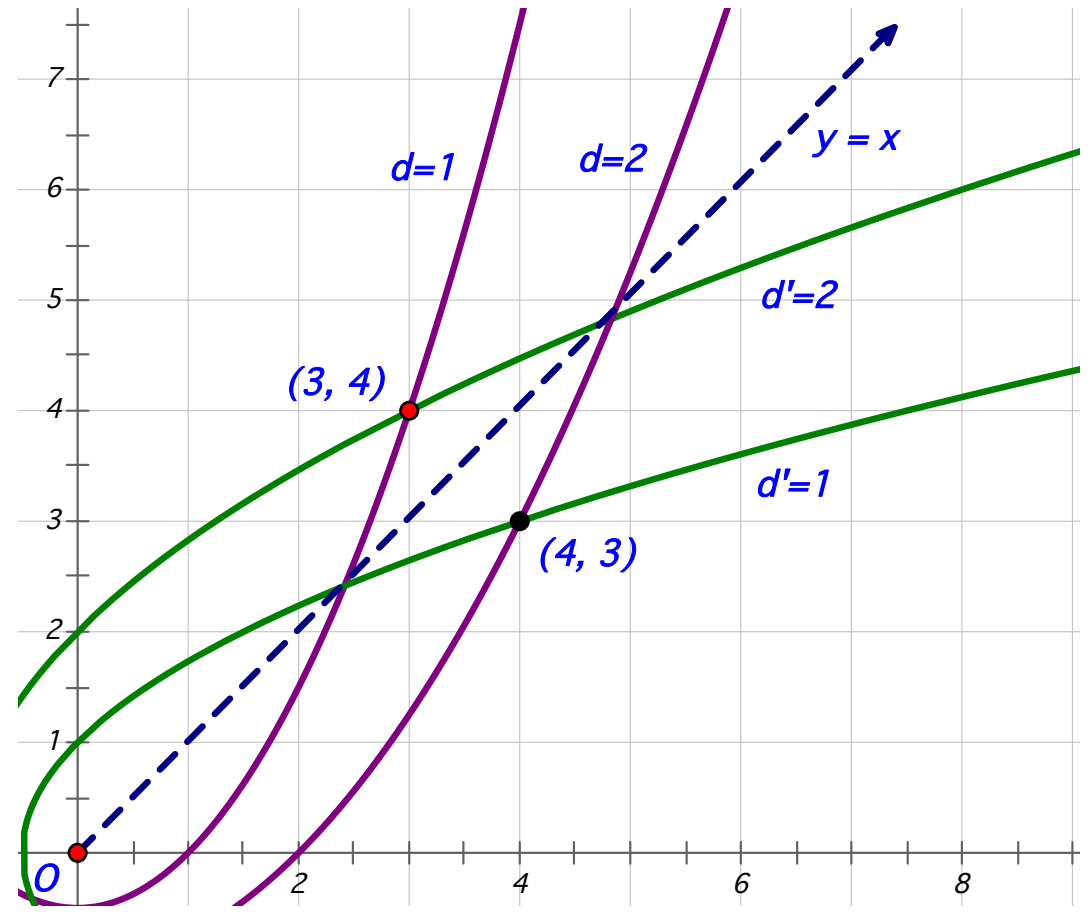

There is another $d$ value associated with the PPT $(a, b, c)$, $a<b<c$, which is determined by the PPT form $(b, a, c)$, namely $d' = c-a$. So the point $(a, b)$ will be at the intersection of two parabolas, one for $d = c-b$, and one for $d' = c-a$. The points $(a, b)$ and $(b, a)$ will each be at the intersection of two parabolas, one for $d'$ and one for $d$. See the graph here, where $d = 1, 2$, and $d' = 1, 2$.

***Theorem 8. Parabolas with d = 1, 2.***

*a.) If $(a, b, c)$ is a PPT with $a<b<c$, and $d = 1$, then $(a, b)$ is a point on the parabola with equation $y = (x^2-1)/2$, denoted by $d = 1$, and $(b, a)$ is a point on the reflected parabola with equation $x = (y^2-1)/2$, which is denoted $d' = 1$.*

*b.) If $(b, a, c)$ is a PPT with $a<b<c$, and $d = 2$, then $(b, a)$ is a point on the parabola with equation $y = (x^2-4)/4$, denoted by $d = 2$, and $(a, b)$ is a point on the reflected parabola with equation $x = (y^2-4)/4$, which is denoted $d' = 2$.*

The graph of all of the parabolas for the allowable values $d, d' = 1, 2, 8, 9, 18, 25, 32, 49$, and $50$ is shown below. This is the $180x180$ lower left corner of the red and black figure above.

This graph shows many reflected points, such as the pairs $(9,40)$ and $(40,9)$. Point $(9,40)$ is at the intersection of parabolas for $d = 1$ and $d' = 32$, and point $(40,9)$ is at the intersection of the parabolas for $d' = 1$ and $d = 32$.

For the point *(9,40)* we have the parabolas $y = (x^2 -1)/2$, where $d=1$, and $x = (y^2 -32^2)/64$, where $d' = 32$, and for the point *(40,9)* we have parabolas $y = (x^2 -32^2)/64$, where $d=32$, and $x = (y^2 -1)/2$, where $d'=1$.

Theorem 8 can be generalized to other values of *d.* In general, given a PPT *(a,b,c)*, there are two parabolas determined. One which contains the point (a, b), denoted by *d*, where $d = c - b$, and one for it's reflected point *(b, a)*, denoted by *d'*, where $d' = c - a$. These parabolas are given by the equations $y = (x^2 -d^2)/2d$, and $x = (y^2 -d'^2)/2d'$, respectively.

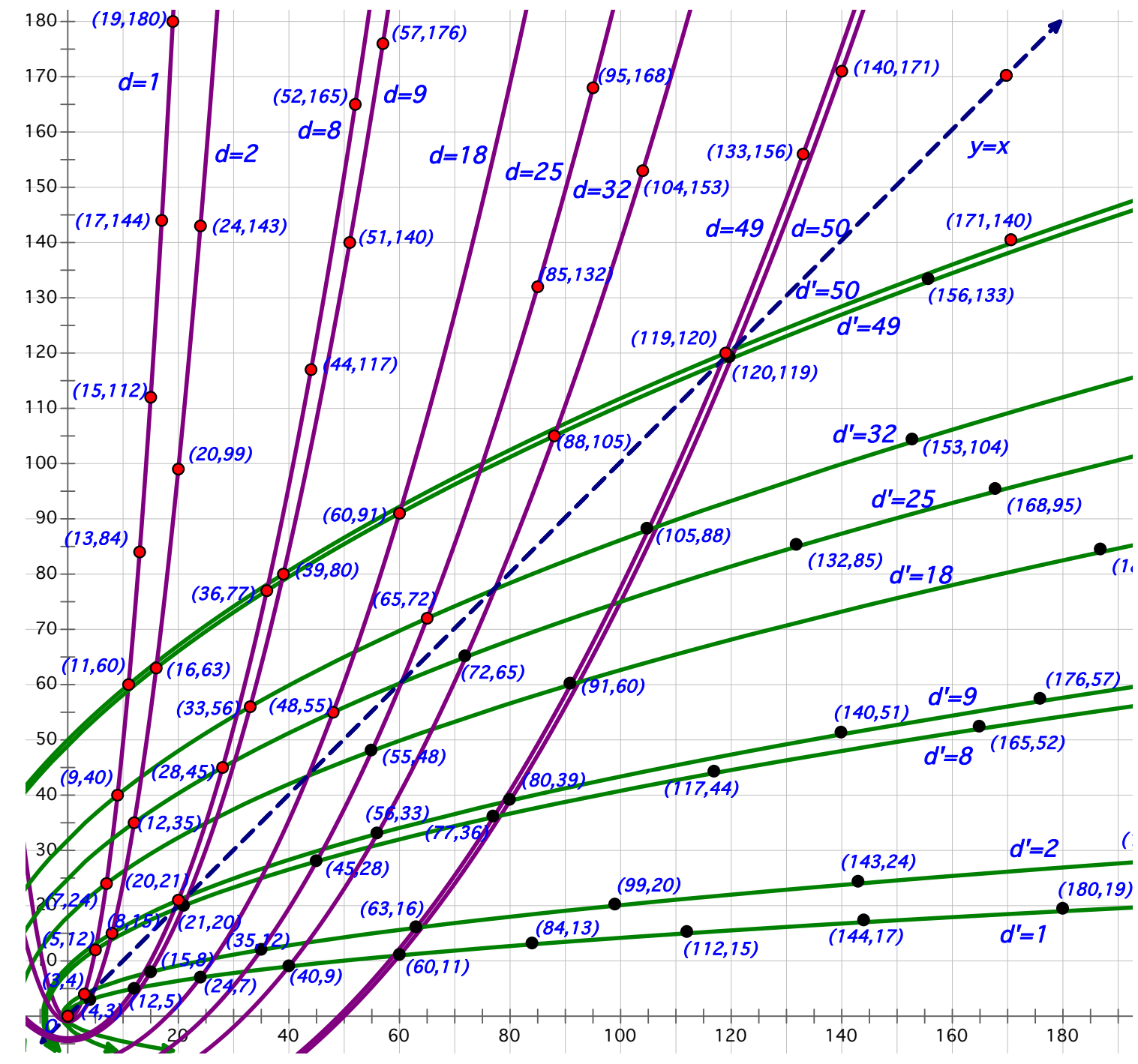

Any integer point *(s, t)* on a parabola for some *d* value will satisfy the equation for that parabola, $y = (x^2 -d^2)/2d$, so $t = (s^2 -d^2)/2d$. Also *(s,t,u)* is a PPT for $d = u - t$, where $u = t + d = (s^2 + d^2)/2d$. Thus each point *(a, b)* in the graph of the PPT *(a,b,c)* is on a parabola for some value of *d,* and also on a parabola for the associated value of *d'.* This explains some of the parabolic patterns seen in the plot of PPTs with legs $a, b < 10,000$.

The amalgamation of the plots of the results for the allowable *d* values using *Mathematica*, determines a result identical to the main plot above of PPTs with legs $a, b < 10,000$ [15].

Also note that all of the *d* parabolas have focus points at the origin *(0,0)*, while the vertex points are at *(0, -d/2)* (show this).

***Theorem 9. Parabolas with Allowable d Values.***

*If (a, b, c) is a PPT with a<b<c, and allowable values d, d', then (a, b) is a point on the parabolas with equations* $y = (x^2 -d^2)/2d$, *and* $x = (y^2 -d'^2)/2d'$, *and (b, a) is a point at the intersection of the reflection of those parabolas, which have equations* $x = (y^2 -d^2)/2d$, *and* $y = (x^2-d'^2)/2d'$.

Some of the other curves visible in the main red and black plot above of the PPTs with legs $a,b < 10,000$ are discussed in [4] and [14].

There are other methods for finding Pythagorean triples, and their properties, see [4]. The range of applications and inter-connections with other mathematical topics is very wide.

***Hints & Answers.***

p. 1. On the proof of the Pythagoras Equation: The area of blue squares $a^2+b^2$ plus the area of the $4$ white triangles $2ab$ equals the area of the large squares $(a+b)^2$, which also equals the area of the red square $c^2$ plus the area of the $4$ white triangles $2ab$. Thus we have $a^2+b^2+2ab = (a+b)^2 = c^2+2ab$, so … .

p. 2. Since we have $\Delta ABD$ is right, by construction, $a^2+b^2 = t^2$, $t = AD$, and since $\Delta ABC$ satisfies $a^2+b^2 = c^2$, by assumption, it follows that $c = t$, so the triangles are congruent by *SSS.* Thus *ACB* is a right angle, and $\Delta ABC$ is a right triangle.

p. 2. "… it is not possible for a PT to have $a = b$", for if it were true, then we would have the equation $2a^2 = c^2$, which implies $\sqrt{2}\cdot a = c$, which means that $c$ is not an integer.

p. 2. Observe $1$ is equal to $1\cdot 1$, and no other integer factors are possible, so the two factors given are impossible, as they are unequal.

p. 2. If $a = 2$, then we would have $4 = (c-b)(c+b)$. The integer factors of $4$ are $1$, $2$ and $4$, $4 = 2\cdot 2 = 1\cdot 4$. So let $(c-b) = 1$, $(c+b) = 4$, and $c = b+1$ (both factors equal to $2$ is impossible). But then $2b = 3$, contradicting the assumption that $b$ is an integer.

p. 3. Proof of Euclid's Formula follows from $(h^2-k^2)^2 + (2hk)^2 = (h^2+k^2)^2$ by substitution, since $(h^2-k^2)^2 = h^4-2h^2k^2+k^4$, $(2hk)^2 = 4h^2k^2$, and $(h^2+k^2)^2 = h^4+2h^2k^2+k^4$.

p. 3. "…if $h$, $k$ are relatively prime… $(h^2 -k^2, 2hk, h^2+k^2)$ is a PPT."

Assume that $h$, $k$ are relatively prime, and $p$ is a factor of both of $h^2\pm k^2$ and of $2hk$. Then $p$ would be a factor of $(h^2-k^2)+(h^2+k^2) = 2h^2$, and thus either $p = 2$ or $p$ is a factor of $h$. But $p = 2$ is impossible, because $p$ is a factor of $h$. If $p$ is a factor of $h$, then $p$ is a factor of $h^2 - (h^2 -k^2) = k^2$, and thus of $k$. But this contradicts the assumption on $h$, $k$, so $(h^2 -k^2, 2hk, h^2+k^2)$ is a PPT.

p. 3. On the non-existence of $h$, $k$ such that $(9 ,12, 15)$ is generated. If $h = 6$, and $k = 1$, then $a = 35$, and if $h = 3$, and $k = 2$ , then $a = 5$.

p. 4. "… $b = (a^2 -1)/2$, which is then a positive integer."

Since $a$ is odd, and $a >1$, $a^2 = (2m+1)^2$, for some $m > 0$, and $a^2 -1 = 4m^2 + 4m$, so it is divisible by $2$ (actually by $4$).

p. 5. "…consecutive integers are always relatively prime"

Numbers which are divisible by $2$ are at least $2$ units apart, viz. $2, 4, 6, 8, \ldots$ ; numbers which are divisible by $3$ are at least $3$ units apart, $3, 6, 9, 12, \ldots$, and so on.

p. 5. "Explain why these three assumptions lead to a contradiction."

The sum $a^2 + b^2$ has the form $4(h^2 + k^2 + h + k) + 2$, so $4(h^2 + k^2 + h + k) + 2 = 4m^2$. Dividing out the common factor $2$, you have $2(h^2 + k^2 + h + k) + 1 = 2m^2$, a contradiction.

p. 5. "Why must $a$ be even?" If $a$ is odd, then $a^2$ is odd, and $a^2 - 4$ is odd, so $(a^2 - 4)/4$ would not be an integer.

p. 7. "…the choices for $d = 8, 9$ determine both PTs and PPTs in each case."

If $d = 8$, then $b = (a^2 - 64)/16$, so $a > 8$, and $a$ must be even. We must check out cases for $a = 0\bmod 8$, $2\bmod 8$, $4\bmod 8$, and $6\bmod 8$. If $a = 0\bmod 8$, then $a = 8k$, $k >1$, $k$ even, and $b = 4(k^2 -1)$, $b+8 = 4(k^2 +1)$. Since all factors are even, we have PTs for all $k > 1$.

For example $(16, 12, 20), (24, 32, 40), (32, 60, 68), \ldots$ . If $a = 2\bmod 8$, or $6\bmod 8$, the term $b$ is not an integer (see examples, this can and should be backed up by algebraic argument). If $a = 4\bmod 8$, then $a = 8k + 4$, an even number, for all $k > 0$. So $b = 4k^2 + 4k$ $-3$, and $b + 8 = 4k^2 + 4k + 5$, both odd. Thus $b$ & $b+8$ do not have any common factors, and we have PPTs for all $k > 0$, in this case. For example $(12, 5, 13)$, $(20, 21, 29)$.

If $d = 9$, then $b = (a^2 - 81)/18$, for $a$ odd, $a > 9$. Checking examples, you will see that we must have $a = 0 \bmod 9$, $3 \bmod 9$, or $6 \bmod 9$. Algebraic arguments will support these choices. If $a = 0 \bmod 9$, then $a = 9k$, $k > 2$, $k$ odd, so $b = 9(k^2 - 1)/2$, and $b+9 = 9(k^2 + 1)/2$. Thus all $3$ factors have a common factor of $9$, so these triples are PTs. If $a = 3 \bmod 9$, then $a = 9k + 3$, $k$ even, and $b = (9k^2 + 6k - 8)/2$, $b + 9 = (9k^2 + 6k + 10)/2$. These $3$ values have no common factors, since $b$ is even, but $b + 9$ is odd, and neither factor is divisible by $3$. Thus these triples are PPTs. Similarly, if $a = 6 \bmod 9$, we have PPTs.

p. 8. "…the percent of prime numbers in a set of odd integers such as $\{1, 3, \ldots, n\}$ approaches $0$ as $n$ increases". What does this mean? This can be illustrated by simply counting the number of odd primes in a set of odd integers. For instance, there are *24* primes in the set $\{1, 3, \ldots, 99\}$ of *50* odd integers, so we have *48%* primes. But in the set $\{1, 3, \ldots, 999\}$, there are *167* primes, which is *33.4%,* in the set $\{1, 3, \ldots, 9999\}$ there are *1229* primes for *24.6%,* and in the set $\{1, 3, \ldots, 99999\}$ there are *9592* primes for *19.2%.* The process is slow, but as the size of the set of odd integers increases, the percent of primes does indeed drop reluctantly towards $0$ (see the *Prime Number Theorem* [13]).

p. 9. "it follows that $m$ and $n$ are relatively prime". Let $e>1$ be a common factor of $m$ and $n$. Then $e$ is a common factor of $a$ and $c$, and then also of $(c - a)$ and $(c + a)$. Thus $e^2$ would divide $c^2$ so $e$ would divide $a$ and $b$, contradiction of assumption on $a$, $b$.

p. 9. "…$c$ always ends in *13…*" . Since $c - 1 = (a^2 - 1)/2 = (25(2k+1)^2 - 1)/2 = \ldots = (100(\ldots) + 24)/2 = 50(\ldots) + 12$, we find that $c = 50(\ldots) + 13$. Note that $k^2 + k = k(k+1)$, so if $k=12$ or $k+1 = 13n$, for some $n$, then $c$ has a factor of $13$.

p.11. "… $u$ and $v$ must both be either even or odd", - if one is even and one odd, then both $u+v$ and $v-u$ would both be odd, and hence not divisible by $2$.

p.12. Both triangles have the same length sides, so they are congruent by *SSS*.

p.15. The vertex points are found by setting $y = 0$, and solving for $x$. The focus points are at $(0, k + 1/(4e))$, where $k = -d/2$, the *y-coordinate* of the vertex point, and $e$ is the coefficient of $x$ in the parabola equation, which in our case is $e = 1/2d$. The value $1/4e = d/2$ is the distance between the vertex point and the focus point.

***References.***